\documentclass[12pt,reqno]{amsart}

\usepackage{amssymb}

\usepackage{eucal}

\input{diagrams}

\font\sss=cmss8

\def\cE{{\mathcal E}}

\def\BA{{\mathbb A}}

\def\BZ{{\mathbb Z}}

\def\sD{\mbox{\sf D}}

\def\sT{\mbox{\sf T}}

\def\ast{{\textstyle *}}

\def\C{\operatorname{C}}

\def\D{\sD}
\def\Dc{\sD^{\operatorname{c}}}
\def\deg{\operatorname{deg}}

\def\Df{\sD^{\operatorname{f}}}
\def\dim{\operatorname{dim}}
\def\Dsmall{\mbox{\sss D}}
\def\dual{\operatorname{D}}

\def\Ext{\operatorname{Ext}}

\def\H{\operatorname{H}}

\def\Hom{\operatorname{Hom}}

\def\inf{\operatorname{inf}}

\def\J{\operatorname{J}}

\def\LTensor{\stackrel{\operatorname{L}}{\otimes}}
\def\opp{\operatorname{op}}

\def\RHom{\operatorname{RHom}}

\def\sup{\operatorname{sup}}

\numberwithin{equation}{part}

\newtheorem{Lemma}{Lemma}[section]
\newtheorem{Theorem}[Lemma]{Theorem}
\newtheorem{Proposition}[Lemma]{Proposition}
\newtheorem{Corollary}[Lemma]{Corollary}

\theoremstyle{definition}
\newtheorem{Definition}[Lemma]{Definition}
\newtheorem{Setup}[Lemma]{Setup}

\newtheorem{Remark}[Lemma]{Remark}

\def\AR{Aus\-lan\-der-Rei\-ten}
\def\compact{small}
\def\DGm{DG mo\-du\-le}

\def\DGrm{DG right-mo\-du\-le}

\def\DGSm{DG $S$-mo\-du\-le}
\def\DGRlm{DG left-$R$-mo\-du\-le}

\def\DGRrm{DG right-$R$-mo\-du\-le}
\def\DGSrm{DG right-$S$-mo\-du\-le}
\def\DGRbm{DG left/right-$R$-mo\-du\-le}
\def\DGSbm{DG left/right-$S$-mo\-du\-le}

\def\DGSolm{DG left-$S^{\opp}$-mo\-du\-le}

\def\DGOmegaSdrm{DG right-$\OmegaSd$-mo\-du\-le}

\def\DGCXkbm{DG left/right-$\CXk$-mo\-du\-le}
\def\OmegaSdnaturalrm{right-$\OmegaSd^{\natural}$-mo\-du\-le}
\def\tors{\sT}  
\def\CSdk{\C^{\ast}(S^d)}
\def\HSdk{\H^{\ast}(S^d)}

\def\OmegaSd{A}  
\def\OmegaSdkres{F}  
\def\OmegaSdkresendo{\cE}  
\def\OmegaSdind{Y}  
\def\CSdkind{Z}  
\def\Space{X}  

\def\Spaceprime{\Space^{\prime}}
\def\CXk{\C^{\ast}(\Space)}
\def\CXprimek{\C^{\ast}(\Spaceprime)}
\def\HXk{\H^{\ast}(\Space)}
\def\HiXk{\H^i(\Space)}

\begin{document}

\title[Auslander-Reiten triangles]
{Auslander-Reiten triangles and quivers over topological spaces}

\author{Peter J\o rgensen}
\address{Danish National Library of Science and Medicine, N\o rre
All\'e 49, 2200 K\o \-ben\-havn N, DK--Denmark}
\email{pej@dnlb.dk, www.geocities.com/popjoerg}


\keywords{\AR\ triangle, quiver, cochain differential graded algebra,
Poincar\'e duality, topological space, sphere}

\subjclass[2000]{55P62, 16E45, 16G70}

\begin{abstract} 

In this paper, \AR\ triangles are introduced into algebraic topology,
and it is proved that their existence characterizes Poincar\'e duality
spaces.

Invariants in the form of quivers are also introduced, and \AR\
triangles and quivers over spheres are computed.  

The quiver over the $d$-dimensional sphere turns out to consist of
$d-1$ components, each isomorphic to $\BZ
\BA_{\hspace{0.15mm}\infty}$.  So quivers are sufficiently sensitive
invariants to tell spheres of different dimension apart.

\end{abstract}

\maketitle

\setcounter{section}{-1}
\section{Introduction}
\label{sec:introduction}

In this paper, two concepts from representation theory are introduced
into algebraic topology: \AR\ triangles and invariants in the form of
qui\-vers (that is, directed graphs).

The highlights are that existence of \AR\ triangles characterizes
Poincar\'e duality spaces (theorem~\ref{thm:Poincare_duality_1}), that
\AR\ triangles and quivers over spheres can be computed
(theorems~\ref{thm:AR_triangles_over_spheres}
and~\ref{thm:quivers_over_spheres}), and that quivers are
sufficiently sensitive invariants to tell spheres of different
dimension apart (corollary~\ref{cor:quivers_tell_spheres_apart}).

After this very short survey, let me describe the paper at a more
leisurely pace.

The idea to use methods from the representation theory of finite
dimensional algebras in algebraic topology comes as follows:

If $k$ is a field and $\Space$ is a simply connected topological space
with $\dim_k \H^{\ast}(\Space;k) < \infty$, then the singular cochain
differential graded algebra $\C^{\ast}(\Space;k)$ is equivalent by a
series of quasi-isomorphisms to a
differential graded algebra $R$ which is finite dimensional over $k$,
by the methods of \cite[proof of thm.\ 3.6]{FHT} and
 \cite[exam.\ 6, p.\ 146]{FHTbook}.

Hence it seems obvious to try to study $R$ and thereby
$\C^{\ast}(\Space;k)$ with methods from the representation theory of
finite dimensional algebras.  A natural place to start is with the
derived category of differential graded modules $\D(R)$ which is the
playing ground for homological algebra over $R$.  Note that by
\cite[thm.\ III.4.2]{KrizMayAst}, the category $\D(R)$ is equivalent
to $\D(\C^{\ast}(\Space;k))$.

A number of concepts present themselves which are used to analyze
the structure of derived categories in representation theory.  I will
concentrate on two important ones: \AR\ triangles and invariants in
the form of quivers.  Their definitions are recalled
in~\ref{def:AR_triangles} and~\ref{def:quivers} below, but let me make
some remarks.

\AR\ triangles were introduced by Happel in \cite{HapDerCat}, and are
certain special triangles among the distinguished triangles in a
triangulated category.  They are the triangulated counterpart to \AR\
sequences which pervade re\-pre\-sen\-ta\-ti\-on theory, see for
instance \cite{ARSbook}.  Not all triangulated categories have \AR\
triangles, but those that do enjoy many advantages.  This is expounded
in Happel's book \cite{Hapbook}, but see also his papers
\cite{HapGor} and \cite{HapDerCat}. 

The quiver of an additive category is an important structural
invariant.  The vertices are certain isomorphism classes in the
category and the arrows are determined by certain morphisms.  One can
think of the quiver as an ``X-ray image'' of the category.  Quivers of
additive categories are used extensively in representation theory; an
example is the so-called \AR\ quiver which is a tremendously useful
tool, see \cite{ARSbook}.

\AR\ triangles and quivers are intimately connected: If a suitable
triangulated category has \AR\ triangles, then they give enough
information to compute the quiver of the category, see
lemma~\ref{lem:AR_triangles_determine_quivers}, and they even give the
quiver the extra structure of so-called stable translation quiver, see
definition~\ref{def:translation_quivers} and
corollary~\ref{cor:AR_triangles_give_translation_quivers}.

Now, one can hope that the tools of \AR\ triangles and quivers will be
as useful in studying the derived category $\D(\C^{\ast}(\Space;k))$
as they are in representation theory.  This paper shows that at least
something can be gained:

Section~\ref{sec:AR_triangles_top} considers \AR\ triangles, and proves
(essentially) that they exist in the category
$\Dc(\C^{\ast}(\Space;k))$ if and only if the topological space
$\Space$ has Poincar\'e duality over $k$
(theorem~\ref{thm:Poincare_duality_1}).  Here
$\Dc(\C^{\ast}(\Space;k))$ is the full subcategory of \compact\
objects of $\D(\C^{\ast}(\Space;k))$ (those where $\Hom(M,-)$ commutes
with set indexed coproducts).

Section~\ref{sec:quivers_top} considers the quiver of
$\Dc(\C^{\ast}(\Space;k))$, and proves that it is a weak homotopy
invariant of $\Space$
(proposition~\ref{pro:quiver_weak_homotopy_invariant}). 

Section~\ref{sec:spheres} applies the theory to spheres, and computes
the \AR\ triangles and the quiver of $\Dc(\C^{\ast}(S^d;k))$ for $d
\geq 2$ when the characteristic of $k$ is zero
(theorems~\ref{thm:AR_triangles_over_spheres}
and~\ref{thm:quivers_over_spheres}).  The quiver consists of $d-1$
components, each isomorphic to $\BZ \BA_{\hspace{0.2mm}\infty}$, and
it is observed that hence, the quiver of $\Dc(\C^{\ast}(S^d;k))$ is a
sufficiently sensitive invariant to tell spheres of different
dimension apart (corollary~\ref{cor:quivers_tell_spheres_apart}).

On the way to these results, the indecomposable objects of the
category $\Dc(\C^{\ast}(S^d;k))$ are determined, and it is proved that
each object is the coproduct of uniquely determined indecomposable
objects (proposition~\ref{pro:indecomposables_over_spheres}).  This
gives a fairly accurate picture of $\Dc(\C^{\ast}(S^d;k))$ which may
be of independent interest.

The initial sections~\ref{sec:AR_triangles}
to~\ref{sec:Gorenstein_DGAs} of the paper are organized as follows:
Sections~\ref{sec:AR_triangles} and~\ref{sec:quivers} briefly recall
\AR\ triangles and quivers, and sections~\ref{sec:derived_categories}
to~\ref{sec:Gorenstein_DGAs} develop the theory of \AR\ triangles over
a general differential graded algebra $R$ which has the advantage of
being typographically lighter than $\C^{\ast}(\Space;k)$, and not
mathematically harder.

\medskip
Let me end the introduction by giving some notation.

Throughout the paper, $k$ denotes a field.

Differential Graded Algebras are abbreviated DGAs, and Differential
Graded modules are abbreviated \DGm s. 

Standard notation is used for triangulated categories and for derived
categories and functors of \DGm s over DGAs.  The suspension functor
is denoted $\Sigma$ and the $i$'th cohomology functor is denoted
$\H^i$.  The notation is cohomological (degrees indexed by
superscripts, differentials of degree $+1$).

Module structures are occasionally emphasized by subscript notation.
So for instance, $M_{R,S}$ indicates that $M$ has compatible
right-structures over $R$ and $S$.

Let $S$ be a DGA over $k$.

$S^{\natural}$ denotes the graded algebra obtained by forgetting the
differential of $S$, and if $M$ is a \DGSm\ then $M^{\natural}$
denotes the graded $S^{\natural}$-module obtained by forgetting the
differential of $M$.

The opposite DGA of $S$ is denoted $S^{\opp}$, and is defined by having
the product $s \stackrel{\opp}{\cdot} t = (-1)^{|s||t|}ts$.  \DGSrm s
are identified with \DGSolm s.

$\Dc(S)$ denotes the full subcategory of the derived category $\D(S)$
which consists of \compact\ $M$'s, that is, $M$'s so that $\Hom(M,-)$
commutes with set indexed coproducts.

$\Df(S)$ denotes the full subcategory of $\D(S)$ which consists of
$M$'s with $\dim_k \H\!M < \infty$.

I write 
\[
  \dual(-) = \Hom_k(-,k).
\]
This duality functor is defined on $k$-vector spaces.  It can also be
viewed as defined on modules, graded modules, or \DGm s, and as such
it interchanges left-modules and right-modules.  The functor
$\dual$ induces a duality of categories
\[
  \begin{diagram}[labelstyle=\scriptstyle]
    \Df(S) && \pile{ \rTo^{\dual} \\ \lTo_{\dual} } && \Df(S^{\opp}). \\
  \end{diagram}
\]
Note that $\dual\!S$ is a \DGSbm, like $S$ itself.

If $\sD$ is a triangulated category and $M$ is an object of $\sD$,
then an object of $\sD$ is said to be finitely built from $M$ if it is in the
smallest triangulated subcategory of $\sD$ which contains $M$ and is
closed under retracts.

\section{\AR\ triangles}
\label{sec:AR_triangles}

Sections~\ref{sec:AR_triangles} and~\ref{sec:quivers} are
introductory.

This section recalls the definition of \AR\ triangles and a few of
their properties from \cite{HapDerCat}.

Let $\sD$ be a $k$-linear triangulated category over the field $k$,
where each $\Hom$ space is finite dimensional over $k$ and each
indecomposable object has local endomorphism ring.

\begin{Definition}
\label{def:AR_triangles}
A distinguished triangle 
\[
  M \longrightarrow N
  \stackrel{\nu}{\longrightarrow} P 
  \stackrel{\pi}{\longrightarrow} 
\]
in $\sD$ is called an {\em \AR\ triangle} in $\sD$ if
\begin{enumerate}

  \item  $M$ and $P$ are indecomposable objects.

  \item  $\pi \not= 0$.

  \item  Each morphism $N^{\prime} \longrightarrow P$ which is not a
         retraction factors through $\nu$.

\end{enumerate}
\end{Definition}

Given an indecomposable object $P$, there may or may not exist an \AR\
triangle as in the definition.  But if there does, then it is
determined up to isomorphism by \cite[prop.\ 3.5(i)]{HapDerCat}.  This
allows the following definition.

\begin{Definition}
\label{def:AR_translate}
Given an indecomposable object $P$ of $\sD$.  Suppose that there is an
\AR\ triangle as in definition~\ref{def:AR_triangles}.  Then $M$ is
called the {\em \AR\ translate} of $P$ and denoted $\tau P$.  The
operation $\tau$ is called the {\em \AR\ translation} of $\sD$.  
\end{Definition}

Note that $\tau P$ is only defined up to isomorphism.  

\begin{Definition}
Suppose that for each indecomposable object $P$ of $\sD$, there exists
an \AR\ triangle as in definition~\ref{def:AR_triangles}.  Then $\sD$
is said to {\em have \AR\ triangles}.
\end{Definition}

\section{Quivers}
\label{sec:quivers}

This section recalls the definition of the quiver of an additive
category and its connection with \AR\ triangles.

Let $\sD$ be an additive category.  

A morphism $M \stackrel{\mu}{\longrightarrow} N$ is called {\em
irreducible} if it is neither a section or a retraction, but any
factorization $\mu = \rho \sigma$ has either $\sigma$ a section or
$\rho$ a retraction.

\begin{Definition}
\label{def:quivers}
The {\em quiver} of $\sD$ has as vertices the isomorphism classes
$[M]$ of indecomposable objects of $\sD$. It has one arrow from $[M]$
to $[N]$ when there is an irreducible morphism $M \longrightarrow N$
and no arrows from $[M]$ to $[N]$ otherwise.
\end{Definition}

Now let $\sD$ be a $k$-linear triangulated category where each $\Hom$
space is finite dimensional over $k$ and each indecomposable
object has local endomorphism ring.

In this case the quiver of $\sD$ and the \AR\ triangles in $\sD$ are
intimately connected.  The following result is immediate from
\cite[prop.\ 3.5]{HapDerCat}.

\begin{Lemma}
\label{lem:AR_triangles_determine_quivers}
Let $M \longrightarrow N \longrightarrow P \longrightarrow$ be an \AR\
triangle in $\sD$.  Suppose that $N \cong \coprod_j N_j$ is a
splitting into indecomposable objects, and let $N^{\prime}$ be some
indecomposable object.  Then the following statements are equivalent.
\begin{enumerate}
  \item  There is an irreducible morphism $M \longrightarrow N^{\prime}$.

  \item  There is an irreducible morphism $N^{\prime} \longrightarrow P$.

  \item  There is a $j$ so that $N^{\prime} \cong N_j$.
\end{enumerate}
\end{Lemma}

So if $\sD$ has \AR\ triangles, then knowledge of the \AR\ triangles
gives full knowledge of the quiver of $\sD$.

Moreover, there is the notion of stable translation quiver.

\begin{Definition}
\label{def:translation_quivers}
A quiver is said to be a {\em stable translation quiver} if it is
equipped with a map $\tau$ called the {\em translation}, which sends
vertices to vertices in a way so that the number of arrows from
$\tau[P]$ to $[N^{\prime}]$ equals the number of arrows from
$[N^{\prime}]$ to $[P]$.
\end{Definition}

Lemma~\ref{lem:AR_triangles_determine_quivers} implies that if $\sD$
has \AR\ triangles, then definition~\ref{def:translation_quivers} is
satisfied with $\tau[P] = [M]$.  Note $[M] = [\tau P]$, where $\tau$
now stands for the \AR\ translation of $\sD$, see
definition~\ref{def:AR_translate}.  Hence the following result.

\begin{Corollary}
\label{cor:AR_triangles_give_translation_quivers}
If $\sD$ has \AR\ triangles, then the quiver of $\sD$ is a stable
translation quiver with translation induced by the \AR\ translation of
$\sD$ via $\tau[P] = [\tau P]$.
\end{Corollary}

\section{Derived categories}
\label{sec:derived_categories}

Sections~\ref{sec:derived_categories}, \ref{sec:AR_triangles_DGA},
and~\ref{sec:Gorenstein_DGAs} develop the theory of \AR\ triangles
over a general DGA denoted $R$.  

This section collects some lemmas on derived categories of
\DGm s.

\begin{Setup}
\label{set:blanket_R}
In the rest of the paper, $R$ denotes a DGA over the field $k$
satisfying:
\begin{enumerate}

  \item  $R$ is a cochain DGA, that is, $R^i = 0$ for $i < 0$.

  \item  $R^0 = k$.

  \item  $R^1 = 0$.

  \item  $\dim_k R < \infty$.

\end{enumerate}
Note that $R/R^{\geq 1} \cong k$ is a \DGRbm.
\end{Setup}

First a general result which holds by \cite[thm.\ 5.3]{KelDer}.

\begin{Lemma}
\label{lem:comp_built_from_R}
Let $S$ be a DGA over $k$.  Then the objects of $\Dc(S)$ are exactly
the ones which are finitely built from ${}_{S}S$.
\end{Lemma}

The rest of this section deals with $R$, the DGA from
setup~\ref{set:blanket_R}.  If $M$ is a \DGRlm, then a semi-free
resolution $F \longrightarrow M$ is called {\em minimal} if the
differential $\partial_F$ takes values in $R^{\geq 1}F$, whence $k
\otimes_R F$ and $\Hom_R(F,k)$ have vanishing differentials.  (See
\cite[chp.\ 6]{FHTbook} for general information on semi-free
resolutions.)  The following result is well known, see
\cite[appendix]{FHT} or
\cite[appendix]{FJiid}. 

\begin{Lemma}
\label{lem:resolution}
Let $M$ be a \DGRlm\ for which $u = \inf \{\,i \,\mid\, \H^i\!M \not=
0 \,\}$ is finite, and for which each $\H^i\!M$ is finite dimensional
over $k$.
\begin{enumerate}

  \item  There is a minimal semi-free resolution $F \longrightarrow M$ which
         has a semi-free filtration with quotients as indicated,
\[
  \begin{diagram}[height=3.5ex,width=2ex]
    & & \!\!\!\!\!\!\!\!\!\!\!\!\!
      \lefteqn{\Sigma^{-u}R^{(\gamma_0)}} 
      & & & & & & & & \!\!\!\!\!\!\!\!\!\!\!\!\!\!\!
      \lefteqn{\Sigma^{-u}R^{(\gamma_1)}}
      & & & & & & & & \!\!\!\!\!\!\!\!\!\!\!\!\!\!\!\!\!\!\!\!\!\!
      \lefteqn{\Sigma^{-u-1}R^{(\gamma_2)}}
    & & & & & & \cdots & & & & \\
    & \ruLine & & \luLine &
    & & & & 
    & \ruLine & & \luLine & 
    & & & & 
    & \ruLine & & \luLine & 
    & & & & & & & & \\
    0 & & \subseteq 
      & & F(0)
      & & \subseteq 
      & & L(1)
      & & \subseteq 
      & & F(1)
      & & \subseteq 
      & & L(2)
      & & \subseteq 
      & & F(2) 
      & & \subseteq
      & & \cdots 
      & & \subseteq 
      & & F, \\
    & & & &
    & \luLine & & \ruLine &
    & & & & 
    & \luLine & & \ruLine & 
    & & & & 
    & & & & & & & & \\
    & & & & & & \!\!\!\!\!\!\!\!\!\!\!\!\!\!\!\!\!\!\!\!\!\!\!
      \lefteqn{\Sigma^{-u-1}R^{(\delta_1)}} 
      & & & & & & & & \!\!\!\!\!\!\!\!\!\!\!\!\!\!\!\!\!\!\!\!\!\!\!
      \lefteqn{\Sigma^{-u-2}R^{(\delta_2)}}
      & & & & & & & & & & \cdots & & & & \\
  \end{diagram}
\]
where superscripts $(\gamma_j)$ and $(\delta_j)$ indicate coproducts.
Here each $\gamma_j$ and each $\delta_j$ is finite, and the induced
map 
\[
  \H^i\!F(m) \longrightarrow \H^i\!F
\] 
is an isomorphism for $i \leq u+m$.  Moreover, I have $\gamma_0
\not= 0$.

  \item  In the construction from {\rm (i)}, I have
\[
  F^{\natural} \cong \coprod_{j \leq -u}
  \Sigma^j(R^{\natural})^{(\beta_j)}, 
\]
where each $\beta_j$ is finite.
\end{enumerate}
\end{Lemma}

\begin{Remark}
\label{rmk:F_and_L_comp}
It follows from lemma~\ref{lem:comp_built_from_R} that each $F(m)$ and
each $L(m)$ in lemma~\ref{lem:resolution}(i) is in $\Dc(R)$, because
each step in the semi-free filtration only adds finitely many
$\Sigma^j R$'s.
\end{Remark}

The following truncation lemma uses that $R^0$ is $k$, and is an
exercise in linear algebra.

\begin{Lemma}
\label{lem:truncations}
\begin{enumerate}

  \item  Let $M$ be a \DGRlm\ for which 
         $u = \inf \{\,i \,\mid\, \H^i\!M \not= 0 \,\}$ is finite.
         Then there exists an injective quasi-isomorphism of \DGRlm s
         $U \longrightarrow M$ with $U^j = 0$ for $j < u$.

  \item  Let $N$ be a \DGRlm\ for which 
         $v = \sup \{\,i \,\mid\, \H^i\!N \not= 0 \,\}$ is finite.
         Then there exists an surjective quasi-isomorphism of \DGRlm s
         $N \longrightarrow V$ with $V^j = 0$ for $j > v$.

\end{enumerate}
\end{Lemma}

The following two lemmas show that $\Df(R)$ and $\Dc(R)$ are
categories of the sort for which \AR\ triangles were defined
in~\ref{def:AR_triangles}. 

\begin{Lemma}
\label{lem:Df_Krull_Schmidt}
\begin{enumerate}

  \item  Let $M$ and $N$ in $\Df(R)$ be given.  Then I have
         $\dim_k \Hom_{\Dsmall(R)}(M,N) < \infty$.

  \item  If $M$ is an indecomposable object of
         $\Df(R)$, then the endomorphism ring $\Hom_{\Dsmall(R)}(M,M)$
         is a local ring.

\end{enumerate}
\end{Lemma}

\begin{proof}
(i):  If $N$ is isomorphic to zero in $\Df(R)$, then part (i) of the
lemma is trivial, so I can suppose that $N$ is not isomorphic to zero.

Let $F \longrightarrow M$ and $G \longrightarrow N$ be
semi-free resolutions chosen according to
lemma~\ref{lem:resolution}(i).  Since I have $\dim_k R < \infty$,
lemma~\ref{lem:resolution}(ii) implies $\dim_k F^j < \infty$ and
$\dim_k G^j < \infty$ for each $j$.

As $N$ is in $\Df(R)$ and is non-isomorphic to zero, the same holds
for $G$, so $u = \inf \{\, i \,|\, \H^i\!G \not= 0 \,\}$ and $v = \sup
\{\, i \,|\, \H^i\!G \not= 0 \,\}$ are finite.  By using both parts of
lemma~\ref{lem:truncations}, I can replace $G$ with a truncation
$G^{\prime}$ so that $G^{\prime}$ is concentrated between degrees $u$
and $v$, and so that $G$ and $G^{\prime}$ are connected by two
quasi-isomorphisms.  As $G^{\prime}$ is a truncation of $G$, I have
$\dim_k G^{\prime j} < \infty$ for each $j$, so altogether
$\dim_k G^{\prime} < \infty$ holds.

But $\dim_k F^j < \infty$ for each $j$ and $\dim_k G^{\prime} <
\infty$ imply 
\[
  \dim_k \Hom_R(F,G^{\prime})^j < \infty
\]
for each $j$, and so
\[
  \Hom_{\Dsmall(R)}(M,N) 
  \cong \H^0(\RHom_R(M,N))
  \cong \H^0(\Hom_R(F,G^{\prime}))
\]
also has $\dim_k \Hom_{\Dsmall(R)}(M,N) < \infty$.

\smallskip
\noindent
(ii): By part (i) and \cite[3.2]{Hapbook}, it is enough to see that
idempotent morphisms in $\Df(R)$ split.  But by \cite[prop.\
3.2]{BokNee} they even do so in $\D(R)$ because $\D(R)$ is a
triangulated category with set indexed coproducts.
\end{proof}

\begin{Lemma}
\label{lem:Dc_in_Df}
There is the inclusion $\Dc(R) \subseteq \Df(R)$.
\end{Lemma}

\begin{proof}
This is clear by lemma~\ref{lem:comp_built_from_R} because ${}_{R}R$
is in $\Df(R)$. 
\end{proof}

Finally, some technicalities.

\begin{Lemma}
\label{lem:Hom_sup}
Let $F$ and $N$ be \DGRlm s with
\[
  F^{\natural} 
    \cong \coprod_{j \leq -u} \Sigma^j(R^{\natural})^{(\beta_j)}
\]
where each $\beta_j$ is finite, and with $N^j = 0$ for $j > v$.  Then
\[
  \sup \{\, i \,|\, \H^i(\Hom_R(F,N)) \not= 0 \,\} \leq -u + v.
\]
\end{Lemma}

\begin{proof}
This follows since
\[
  \Hom_R(F,N)^{\natural} 
  \cong \Hom_{R^{\natural}}(F^{\natural},N^{\natural})
  \cong \prod_{j \leq -u}\Sigma^{-j}(N^{\natural})^{(\beta_j)}
\]
is zero in degrees $> -u + v$, because the highest degree contribution
to the product comes from $\Sigma^u(N^{\natural})^{(\beta_{-u})}$
which is certainly zero in degrees $> -u + v$.
\end{proof}

\begin{Lemma}
\label{lem:RHom_sup}
Let $M$ and $N$ be in $\Df(R)$.  Then
\begin{eqnarray*}
  \lefteqn{\sup \{\, i \,|\, \H^i(\RHom_R(M,N)) \not= 0 \,\}}
  \hspace{10ex} & & \\
  & & = - \inf \{\, i \,|\, \H^i\!M \not= 0 \,\}
        + \sup \{\, i \,|\, \H^i\!N \not= 0 \,\}.
\end{eqnarray*}
\end{Lemma}

\begin{proof}
If $M$ or $N$ is isomorphic to zero in $\Df(R)$, then the equation
just says $-\infty = -\infty$, so I can suppose that neither $M$ or
$N$ is isomorphic to zero.  Then $u = \inf \{\, i \,|\, \H^i\!M \not=
0 \,\}$ and $v = \sup \{\, i \,|\, \H^i\!N \not= 0 \,\}$ are finite.

By lemma~\ref{lem:resolution}(i), pick a semi-free resolution $F
\longrightarrow M$ with
\[
  F^{\natural} 
    \cong \coprod_{j \leq -u} \Sigma^j(R^{\natural})^{(\beta_j)}.
\]
By lemma~\ref{lem:truncations}(ii), replace $N$ with a
quasi-isomorphic truncation with $N^j = 0$ for $j > v$.

Since $\RHom_R(M,N) \cong \Hom_R(F,N)$ holds, what I must prove is
\[
  \sup \{\, i \,|\, \H^i(\Hom_R(F,N)) \not= 0 \,\} = -u + v.
\]
Here $\leq$ follows from lemma~\ref{lem:Hom_sup}, so it remains to
show 
\begin{equation}
\label{equ:non_vanishing_H}
  \H^{-u+v}(\Hom_R(F,N)) \not= 0.
\end{equation}

For this, note that the semi-free filtration of $F$
in lemma~\ref{lem:resolution}(i) gives that there is a semi-split
exact sequence of \DGRlm s,
\begin{equation}
\label{equ:semi_split_sequence}
  0 
  \rightarrow \Sigma^{-u}R^{(\gamma_0)}
  \longrightarrow F
  \longrightarrow F^{\prime}
  \rightarrow 0,
\end{equation}
with $\gamma_0 \not= 0$.  Here the left hand term is just $F(0)$, and
$F^{\prime}$ is the quotient $F / F(0)$.  From the part of the
semi-free filtration which continues up from $F(0)$ follows that
$F^{\prime}$ is graded free with
\begin{equation}
\label{equ:Fprime_splitting}
  (F^{\prime})^{\natural} 
    \cong \coprod_{j \leq -u} \Sigma^j(R^{\natural})^{(\beta^{\prime}_j)}.
\end{equation}

Since the sequence \eqref{equ:semi_split_sequence} is semi-split,
applying $\Hom_R(-,N)$ gives a short exact sequence of complexes,
\[
  0
  \rightarrow \Hom_R(F^{\prime},N)
  \longrightarrow \Hom_R(F,N)
  \longrightarrow \Sigma^u N^{(\gamma_0)}
  \rightarrow 0.
\]
The long exact cohomology sequence of this contains
\[
  {\scriptstyle \H^{-u+v}(\Hom_R(F,N)) }
  \longrightarrow {\scriptstyle \H^{-u+v}(\Sigma^u N^{(\gamma_0)}) }
  \longrightarrow {\scriptstyle \H^{-u+v+1}(\Hom_R(F^{\prime},N)). }
\]
The middle term is $\H^v(N^{(\gamma_0)})$ which is non-zero.  The last
term is zero because lemma~\ref{lem:Hom_sup} and equation
\eqref{equ:Fprime_splitting} imply
\[
  \sup \{\, i \,|\, \H^i(\Hom_R(F^{\prime},N)) \not= 0 \,\} 
  \leq -u + v.
\]
But then the first term is non-zero, proving equation
\eqref{equ:non_vanishing_H}.
\end{proof}

\begin{Lemma}
\label{lem:approximation_by_comp}
Let $N^{\prime}$ in $\Df(R)$ and $v$ in $\BZ$ be given.  Then there exists a
distinguished triangle in $\D(R)$,
\[
  F \longrightarrow N^{\prime} \longrightarrow Q \longrightarrow,
\]
so that $F$ is in $\Dc(R)$ and so that $Q$ is in $\Df(R)$ with $\inf
\{\, i \,|\, \H^i\!Q \not= 0 \,\} \geq v$.
\end{Lemma}

\begin{proof}
Let me use lemma~\ref{lem:resolution}(i) to pick a semi-free resoution
$G \longrightarrow N^{\prime}$.  The semi-free filtration
in lemma~\ref{lem:resolution}(i) gives $G(m) \longrightarrow G$ with
$\H^i\!G(m) \longrightarrow \H^i\!G$ an isomorphism for $i \leq u +
m$, and with $G(m)$ in $\Dc(R)$ by remark~\ref{rmk:F_and_L_comp}.  By
picking $m$ large enough, I can arrange that $\H^i\!G(m)
\longrightarrow \H^i\!G$ is an isomorphism for $i \leq v$.

But then the composition $G(m) \longrightarrow G \longrightarrow
N^{\prime}$ also has $\H^i\!G(m) \longrightarrow \H^i\!N^{\prime}$ an
isomorphism for $i \leq v$, and completing to a distinguished triangle
\begin{equation}
\label{equ:approximation_by_comp_triangle}
  G(m) \longrightarrow N^{\prime} \longrightarrow Q \longrightarrow, 
\end{equation}
the long exact cohomology sequence proves 
$\inf \{\, i \,|\, \H^i\!Q \not= 0 \,\} \geq v$.

So \eqref{equ:approximation_by_comp_triangle} can be used
as the lemma's $F \longrightarrow N^{\prime} \longrightarrow Q
\longrightarrow$.
\end{proof}

\section{\AR\ triangles over a DGA}
\label{sec:AR_triangles_DGA}

Recall $R$, the DGA from setup~\ref{set:blanket_R}.  This section
gives a criterion for the existence of \AR\ triangles in $\Dc(R)$
(proposition~\ref{pro:AR_triangles_existence}), and a formula for
\AR\ triangles when they exist
(proposition~\ref{pro:AR_triangles_computation}). 

Note that by lemmas~\ref{lem:Df_Krull_Schmidt} and~\ref{lem:Dc_in_Df},
both $\Df(R)$ and $\Dc(R)$ are categories of the sort for which
\AR\ triangles were defined in~\ref{def:AR_triangles}, so the concept
makes sense for them.

\begin{Lemma}
\label{lem:AR_triangles}
Let $P$ be an indecomposable object of $\Dc(R)$.  Then there is an
\AR\ triangle in $\Df(R)$,
\[
  \Sigma^{-1}(\dual\!R \LTensor_R P) 
  \longrightarrow N 
  \longrightarrow P
  \longrightarrow.
\]
\end{Lemma}

\begin{proof}
This is a consequence of the theory developed in \cite{KrauAR}:

The natural equivalence
\begin{equation}
\label{equ:pre_Krause}
  \dual(\Hom_{\Dsmall(R)}(P,-)) 
  \simeq \Hom_{\Dsmall(R)}(-,\dual\!R \LTensor_R P)
\end{equation}
holds for $P$ equal to ${}_{R}R$, and therefore also holds for the
given $P$ because $P$ is in $\Dc(R)$ and therefore finitely built from
${}_{R}R$ by lemma~\ref{lem:comp_built_from_R}.  
Let
\[
  \Gamma = \Hom_{\Dsmall(R)}(P,P)
\]
be the endomorphism ring of $P$, and rewrite the left hand side of
\eqref{equ:pre_Krause} to
\[
  \Hom_{\Gamma^{\opp}}(\Hom_{\Dsmall(R)}(P,-),\dual\!\Gamma).
\]
This gives the natural equivalence
\begin{equation}
\label{equ:Krause}
  \Hom_{\Gamma^{\opp}}(\Hom_{\Dsmall(R)}(P,-),\dual\!\Gamma) 
  \simeq \Hom_{\Dsmall(R)}(-,\dual\!R \LTensor_R P).
\end{equation}

Now, since $P$ is an indecomposable object of $\Dc(R)$ and hence of
$\Df(R)$, the endomorphism ring $\Gamma$ is finite dimensional over
$k$ and local by lemma~\ref{lem:Df_Krull_Schmidt}.  The unique simple
left-$\Gamma$-module ${}_{\Gamma}S = {}_{\Gamma}(\Gamma/\J(\Gamma))$
is also finite dimensional over $k$ and has only trivial submodules.
Hence the dual module $(\dual\!S)_{\Gamma}$ has only trivial quotient
modules, so must be the unique simple right-$\Gamma$-module,
$(\Gamma/\J(\Gamma))_{\Gamma}$.  Moreover, the projective cover
${}_{\Gamma}\Gamma \longrightarrow {}_{\Gamma}S$ dualizes to an
injective envelope $(\dual\!S)_{\Gamma}
\longrightarrow (\dual\!\Gamma)_{\Gamma}$.  

So $(\dual\!\Gamma)_{\Gamma}$ is the injective envelope of the unique
simple right-$\Gamma$-module $(\Gamma/\J(\Gamma))_{\Gamma}$.
Therefore, by \cite[def.\ 2.1, thm.\ 2.2, and lem.\ 2.3]{KrauAR}, the
equivalence \eqref{equ:Krause} implies that there is a distinguished
triangle in $\D(R)$,
\begin{equation}
\label{equ:candidate_AR_triangle}
  \Sigma^{-1}(\dual\!R \LTensor_R P) 
  \longrightarrow N
  \stackrel{\nu}{\longrightarrow} P
  \stackrel{\pi}{\longrightarrow},
\end{equation}
satisfying, among other things,
\begin{enumerate}

  \item  $\Sigma^{-1}(\dual\!R \LTensor_R P)$ is an indecomposable
         object of $\D(R)$  (as is $P$ by assumption).

  \item  $\pi \not= 0$.

  \item  Each morphism $N^{\prime} \longrightarrow P$ in $\D(R)$ which
         is not a retraction factors through $\nu$.

\end{enumerate}
(In \cite{KrauAR}, the triangle \eqref{equ:candidate_AR_triangle} is
called an \AR\ triangle, but his definition of this concept differs
from mine.) 

Moreover, \eqref{equ:candidate_AR_triangle} is in $\Df(R)$: As $P$ is
finitely built from ${}_{R}R$, it follows that $\dual\!R
\LTensor_R P$ is finitely built from ${}_{R}(\dual\!R)$.  
But then $\dual\!R \LTensor_R P$ is in $\Df(R)$ because
${}_{R}(\dual\!R)$ is in $\Df(R)$.  And $P$ is also in $\Df(R)$ by
lemma~\ref{lem:Dc_in_Df}.  So both end terms in 
\eqref{equ:candidate_AR_triangle} are in $\Df(R)$, and the
long exact cohomology sequence then proves the same for the middle term.

Together, these properties of the distinguished triangle
\eqref{equ:candidate_AR_triangle} imply that it is an
\AR\ triangle in $\Df(R)$; cf.\ definition~\ref{def:AR_triangles}.
\end{proof}

\begin{Lemma}
\label{lem:AR_triangles_extend}
If
\begin{equation}
\label{equ:AR_triangles_extend}
  M \longrightarrow N 
  \stackrel{\nu}{\longrightarrow} P
  \stackrel{\pi}{\longrightarrow}
\end{equation}
is an \AR\ triangle in $\Dc(R)$, then it is also one in $\Df(R)$.
\end{Lemma}

\begin{proof}
When viewed in $\Dc(R)$, the \AR\ triangle
\eqref{equ:AR_triangles_extend} is characterized by satisfying
conditions (i) to (iii) of definition~\ref{def:AR_triangles}.
Clear\-ly, when viewed in $\Df(R)$, the triangle again satisfies
conditions (i) and (ii).  It remains to check condition (iii).

So suppose that 
\[
  N^{\prime} \stackrel{\nu^{\prime}}{\longrightarrow} P 
\]
is a non-retraction in $\Df(R)$.  I must show that $\nu^{\prime}$
factors through $\nu$, which is equivalent to 
\begin{equation}
\label{equ:goal}
  \pi \nu^{\prime} = 0.
\end{equation}

To prove this, let me first write $v = \sup \{\, i \,|\, \H^i\!M \not=
0 \,\}$.  This is finite because $M$ is indecomposable in $\Dc(R)$,
hence not isomorphic to zero.  By
lemma~\ref{lem:approximation_by_comp} there is a distinguished
triangle in $\D(R)$,
\begin{equation}
\label{equ:approximation_triangle}
  F
  \stackrel{\varphi}{\longrightarrow} N^{\prime} 
  \longrightarrow Q
  \longrightarrow,
\end{equation}
with $F$ in $\Dc(R)$ and $Q$ in $\Df(R)$ with 
\begin{equation}
\label{equ:Q_inf_condition}
  \inf \{\, i \,|\, \H^i\!Q \not= 0 \,\} \geq v.
\end{equation}

Here I claim
\begin{equation}
\label{equ:part_goal}
  \pi\nu^{\prime}\varphi = 0, 
\end{equation}
which is a first approximation to equation \eqref{equ:goal}.  To see
this, note that as $F$ is in $\Dc(R)$ and as
\eqref{equ:AR_triangles_extend} is an \AR\ triangle in $\Dc(R)$, it is
enough to see that $\nu^{\prime} \varphi$ is not a retraction.  But it
is not for if there were a section $P
\stackrel{\sigma}{\longrightarrow} F$ with $(\nu^{\prime} \varphi)
\sigma = 1_P$, then $\nu^{\prime} (\varphi \sigma) = 1_P$ would mean
that $\nu^{\prime}$ had the section $\varphi \sigma$, but
$\nu^{\prime}$ is not a retraction.

Next, note $\sup \{\, i \,|\, \H^i(\Sigma M) \not= 0 \,\} =
\sup \{\, i \,|\, \H^i\!M \not= 0 \,\} - 1 = v - 1$.  
Using this and equation \eqref{equ:Q_inf_condition} proves $\leq$ in
\begin{eqnarray*}
  \lefteqn{ \sup \{\, i \,|\, \H^i(\RHom_R(Q,\Sigma M)) \not= 0 \,\} }
    \hspace{10ex} & & \\
  & & \stackrel{\rm (a)}{=}
        - \inf \{\, i \,|\, \H^i\!Q \not= 0 \,\}
        + \sup \{\, i \,|\, \H^i(\Sigma M) \not= 0 \,\} \\
  & & \leq - v + v - 1 \\
  & & = -1,
\end{eqnarray*}
where (a) is by lemma~\ref{lem:RHom_sup}.  Hence the = in
\begin{equation}
\label{equ:vanishing_Hom2}
  \Hom_{\Dsmall(R)}(Q,\Sigma M) \cong \H^0(\RHom_R(Q,\Sigma M)) = 0.
\end{equation}

However, the distinguished triangle \eqref{equ:approximation_triangle}
gives a long exact sequence containing
\[
  \Hom_{\Dsmall(R)}(Q,\Sigma M)
  \longrightarrow \Hom_{\Dsmall(R)}(N^{\prime},\Sigma M)
  \longrightarrow \Hom_{\Dsmall(R)}(F,\Sigma M),
\]
where $\pi\nu^{\prime}$ is an element in the middle term.  The right
hand map sends $\pi\nu^{\prime}$ to $\pi\nu^{\prime}\varphi$ which is
zero by equation \eqref{equ:part_goal}.  So $\pi\nu^{\prime}$ is in
the image of the left hand map, and this image is zero by equation
\eqref{equ:vanishing_Hom2}.  This proves equation \eqref{equ:goal}. 
\end{proof}

\begin{Proposition}
\label{pro:AR_triangles_existence}
The category $\Dc(R)$ has \AR\ triangles if and only if
${}_{R}(\dual\!R)$ is in $\Dc(R)$.
\end{Proposition}

\begin{proof}
On one hand, suppose that ${}_{R}(\dual\!R)$ is in $\Dc(R)$.  Let $P$
be an indecomposable object of $\Dc(R)$.  Then
lemma~\ref{lem:AR_triangles} gives an \AR\ triangle in $\Df(R)$.  In
the present situation, I claim that the triangle is in fact in
$\Dc(R)$, from which follows readily that it is an
\AR\ triangle in $\Dc(R)$; cf.\ definition~\ref{def:AR_triangles}.

To see this, note that as $P$ is in $\Dc(R)$, it is finitely built
from ${}_{R}R$ by lemma~\ref{lem:comp_built_from_R} whence $\dual\!R
\LTensor_R P$ is finitely built from ${}_{R}(\dual\!R)$.  But since
${}_{R}(\dual\!R)$ is in $\Dc(R)$, it is also finitely built from
${}_{R}R$.  All in all, $\dual\!R \LTensor_R P$ is finitely built from
${}_{R}R$, so is in $\Dc(R)$.  But as both $\dual\!R \LTensor_R P$ and $P$
are in $\Dc(R)$, so is the middle term in the distinguished triangle
from lemma~\ref{lem:AR_triangles}, so the triangle is in $\Dc(R)$.

On the other hand, suppose that $\Dc(R)$ has \AR\ triangles.  Let
${}_{R}R \cong \coprod_j R_j$ be a splitting into indecomposable
objects of $\Dc(R)$; such a splitting clearly exists since $\dim_k
\H\!R < \infty$.  Now, for each $j$ there is an \AR\ triangle in
$\Dc(R)$, 
\[
  M_j \longrightarrow N_j \longrightarrow R_j \longrightarrow,
\]
and by lemma~\ref{lem:AR_triangles_extend} this is even an \AR\
triangle in $\Df(R)$.  Also for each $j$ there is an \AR\ triangle in
$\Df(R)$, 
\[
  \dual\!R \LTensor_R R_j 
  \longrightarrow N^{\prime}_j
  \longrightarrow R_j
  \longrightarrow,
\]
by lemma~\ref{lem:AR_triangles}.  

However, the two \AR\ triangles have the same right hand end term,
$R_j$, so by \cite[prop.\ 3.5(i)]{HapDerCat} they are isomorphic.  In
particular, the left hand end terms are isomorphic, so
$M_j \cong \dual\!R \LTensor_R R_j$.  Hence
\[
  \coprod_j M_j \cong \coprod_j \dual\!R \LTensor_R R_j
  \cong \dual\!R \LTensor_R \coprod_j R_j 
  \cong \dual\!R \LTensor_R R \cong {}_{R}(\dual\!R),
\]
and here the left hand side is in $\Dc(R)$ so ${}_{R}(\dual\!R)$
is also in $\Dc(R)$.
\end{proof}

The following result complements lemma~\ref{lem:AR_triangles}.

\begin{Proposition}
\label{pro:AR_triangles_computation}
Suppose that $\Dc(R)$ has \AR\ triangles.
\begin{enumerate}

  \item  Let $P$ be an indecomposable object of $\Dc(R)$.  Then there
         is an \AR\ triangle in $\Dc(R)$,
\[
  \Sigma^{-1}(\dual\!R \LTensor_R P) 
  \longrightarrow N 
  \longrightarrow P
  \longrightarrow.
\]

  \item  The \AR\ translation of $\Dc(R)$ is given by 
\[
  \tau(-) = \Sigma^{-1}(\dual\!R \LTensor_R -).
\]

\end{enumerate}
\end{Proposition}

\begin{proof}
(i): The distinguished triangle here is the one from
lemma~\ref{lem:AR_triangles}, so is an \AR\ triangle in $\Df(R)$.  The
first part of the proof of
proposition~\ref{pro:AR_triangles_existence} shows that it is also an
\AR\ triangle in $\Dc(R)$ provided ${}_{R}(\dual\!R)$ is in $\Dc(R)$.
And this holds by proposition~\ref{pro:AR_triangles_existence} because
$\Dc(R)$ has \AR\ triangles.

\smallskip
\noindent
(ii):  This is immediate from part (i); cf.\
definition~\ref{def:AR_translate}. 
\end{proof}

\section{Poincar\'e duality DGAs}
\label{sec:Gorenstein_DGAs}

Recall $R$, the DGA from setup~\ref{set:blanket_R}.  This section
considers the situation where ${}_{R}(\dual\!R)$ is in $\Dc(R)$ and
$(\dual\!R)_{R}$ is in $\Dc(R^{\opp})$, cf.\
proposition~\ref{pro:AR_triangles_existence}.
Theorem~\ref{thm:Gorenstein_DGAs} shows that this is equivalent to
$\H\!R$ having Poincar\'e duality.

Note that by the proof of theorem~\ref{thm:Gorenstein_DGAs}, it is
also equivalent to $R$ being a so-called Gorenstein DGA; cf.\
\cite{FHT}. 

\begin{Theorem}
\label{thm:Gorenstein_DGAs}
With $d = \sup \{\, i \,|\, \H^i\!R \not= 0 \,\}$,
the following conditions are equivalent.
\begin{enumerate}

  \item  ${}_{R}(\dual\!R)$ is in $\Dc(R)$ 
         and $(\dual\!R)_{R}$ is in $\Dc(R^{\opp})$. 

  \item  There are isomorphisms of graded $\H\!R$-modules
         ${}_{\H\!R}(\dual\!\H\!R) \cong {}_{\H\!R}(\Sigma^d \H\!R)$
         and $(\dual\!\H\!R)_{\H\!R} \cong (\Sigma^d \H\!R)_{\H\!R}$.

\end{enumerate}
\end{Theorem}

\begin{proof}
To facilitate the proof, here are three more conditions each of which
is equivalent to the ones in the theorem.
\begin{enumerate}
\setcounter{enumi}{2}

  \item  $\dim_k \Ext_R(k,R) < \infty$ 
         and $\dim_k \Ext_{R^{\opp}}(k,R) < \infty$.

  \item  There are isomorphisms of graded $k$-vector spaces
         $\Ext_R(k,R) \cong \Sigma^{-d}k$ and $\Ext_{R^{\opp}}(k,R)
         \cong \Sigma^{-d}k$.

  \item  There are isomorphisms ${}_{R}(\dual\!R) \cong
         {}_{R}(\Sigma^d R)$ in $\D(R)$ and $(\dual\!R)_R \cong
         (\Sigma^d R)_R$ in $\D(R^{\opp})$.

\end{enumerate}

\smallskip
\noindent
(i) $\Rightarrow$ (iii):  Duality gives
\begin{equation}
\label{equ:Ext_of_duals}
  \Ext_{R^{\opp}}(k,R) \cong \Ext_R(\dual\!R,\dual\!k) 
  \cong \Ext_R(\dual\!R,k) = (*).
\end{equation}
When (i) holds, lemma~\ref{lem:comp_built_from_R} implies that
${}_{R}(\dual\!R)$ is finitely built from ${}_{R}R$, and then
$\Ext_R(\dual\!R,k)$ is finite dimensional over $k$ since
$\Ext_R(R,k) \cong k$ is finite dimensional over $k$.  Equation
\eqref{equ:Ext_of_duals} then shows that $\Ext_{R^{\opp}}(k,R)$ is
finite dimensional over $k$.  This gives half of (iii), and the other
half follows by symmetry.

\smallskip
\noindent
(iii) $\Rightarrow$ (i): Let $F \longrightarrow {}_{R}(\dual\!R)$ be a
minimal semi-free resolution picked according to
lemma~\ref{lem:resolution}(i).  Continuing the computation from
equation \eqref{equ:Ext_of_duals} gives
\begin{equation}
\label{equ:Ext_of_duals_2}
  (*) = \H(\RHom_R(\dual\!R,k))
  \cong \H(\Hom_R(F,k)) 
  \cong \Hom_{R^{\natural}}(F^{\natural},k^{\natural}),
\end{equation}
where the last $\cong$ is by minimality of $F$.  When (iii) holds,
$\Ext_{R^{\opp}}(k,R)$ is finite dimensional over $k$, and equations
\eqref{equ:Ext_of_duals} and \eqref{equ:Ext_of_duals_2} then show that
$\Hom_{R^{\natural}}(F^{\natural},k^{\natural})$ is finite dimensional
over $k$.  This means that there are only finitely many summands
$\Sigma^j R^{\natural}$ in $F^{\natural}$, so the semi-free filtration
of $F$ in lemma~\ref{lem:resolution}(i) must terminate after finitely many
steps.  So $F$ and therefore ${}_{R}(\dual\!R)$ is finitely built from
${}_{R}R$, whence ${}_{R}(\dual\!R)$ is in $\Dc(R)$.  This gives half
of (i), and the other half fol\-lows by symmetry.

\smallskip
\noindent
(iii) $\Rightarrow$ (iv):  Assume (iii).  The proof that (iii) implies
(i) considered a minimal semi-free resolution $F
\longrightarrow {}_{R}(\dual\!R)$ obtained from
lemma~\ref{lem:resolution}(i), and proved that the
semi-free filtration of $F$ in~\ref{lem:resolution}(i) terminates
after finitely many steps.  But then $F$ must be bounded because
$\dim_k R < \infty$ implies that $R$ itself is bounded.  
Now, the dual of $F \longrightarrow {}_{R}(\dual\!R)$ is
\[
  R_R \cong \dual({}_{R}(\dual\!R)) \longrightarrow \dual\!F,
\]
and this is a $K$-injective resolution of $R_R$ where $\dual\!F$ is
bounded because $F$ is.

Also, lemma~\ref{lem:resolution} gives that ${}_{R}k$ has a
semi-free resolution $G \longrightarrow {}_{R}k$ with 
\begin{equation}
\label{equ:semi_free_resolution_of_k}
  G^{\natural}
    \cong \coprod_{j \leq 0} \Sigma^j(R^{\natural})^{(\beta_j)}
\end{equation}
and each $\beta_j$ finite.

The existence of these resolutions implies that the canonical morphism
\begin{equation}
\label{equ:canonical_isomorphism}
  {}_{R}k \longrightarrow \RHom_{R^{\opp}}(\RHom_R(k,R),R)
\end{equation}
is an isomorphism by \cite[sec.\ 1, thm.\ 1]{Apas}.  Hence
\begin{eqnarray*}
  0 & = & \sup \{\, i \,|\, \H^i({}_{R}k) \not= 0 \,\} \\
    & = & \sup \{\, i \,|\, 
          \H^i(\RHom_{R^{\opp}}(\RHom_R(k,R),R)) \not= 0 \,\} \\
    & \stackrel{\rm (a)}{=} & 
          - \inf \{\, i \,|\, \H^i(\RHom_R(k,R)) \not= 0 \,\}
          + \sup \{\, i \,|\, \H^i\!R \not= 0 \,\} \\
    & = & - \inf \{\, i \,|\, \H^i(\RHom_R(k,R)) \not= 0 \,\} + d,
\end{eqnarray*}
where (a) follows from lemma~\ref{lem:RHom_sup}.  The lemma can be used
because (iii) implies that $\RHom_R(k,R)$ is in $\Df(R^{\opp})$,
while $R_R$ is certainly in $\Df(R^{\opp})$.  This shows
\[
  \inf \{\, i \,|\, \H^i(\RHom_R(k,R)) \not= 0 \,\} = d.
\]

On the other hand, 
\begin{eqnarray*}
  \lefteqn{\sup \{\, i \,|\, \H^i(\RHom_R(k,R)) \not= 0 \,\}} \;\;\;\; & & \\
  & \stackrel{\rm (b)}{=} & 
          - \inf \{\, i \,|\, \H^i\!k \not= 0 \,\}
          + \sup \{\, i \,|\, \H^i\!R \not= 0 \,\} \\
  & = & d,
\end{eqnarray*}
where (b) is again by lemma~\ref{lem:RHom_sup}.

The last two equations show that $\H(\RHom_R(k,R))$ is concentrated in
degree $d$.  Lemma~\ref{lem:truncations} now implies that
$\RHom_R(k,R)$ itself is isomorphic in $\D(R^{\opp})$ to a \DGRrm\
concentrated in degree $d$.  This \DGRrm\ must have the form
$\Sigma^{-d}k_R^{(\alpha)}$, so I get
\[
  \RHom_R(k,R) \cong \Sigma^{-d}k_R^{(\alpha)}.
\]
Inserting this into equation \eqref{equ:canonical_isomorphism} proves
$\alpha = 1$, so all in all 
\[
  \RHom_R(k,R) \cong \Sigma^{-d}k_R
\]
holds.  Taking cohomology gives half of (iv).  The other half follows
by symmetry.

\smallskip
\noindent
(iv) $\Rightarrow$ (iii):  This is clear.

So now, the equivalence of (i), (iii), and (iv) is established.  I
close the proof by establishing the equivalence of (ii), (iv), and
(v). 

\smallskip
\noindent
(ii) $\Rightarrow$ (iv):  This is immediate from the Eilenberg-Moore
spectral sequence
\[
  E_2^{pq} = \Ext_{\H\!R}^p(k,\H\!R)^q 
  \Rightarrow \Ext_R^{p+q}(k,R)
\]
as found in \cite[1.3(2)]{FHT}, and the corresponding spectral
sequence over $R^{\opp}$. 

\smallskip
\noindent
(iv) $\Rightarrow$ (v):  Equation \eqref{equ:Ext_of_duals} gives that
(iv) implies 
\[
  \Ext_R(\dual\!R,k) \cong \Sigma^{-d}k.
\]
Using a minimal semi-free resolution of ${}_{R}(\dual\!R)$, it is easy
to see that this implies half of (v), and the other half
follows by symmetry.

\smallskip
\noindent
(v) $\Rightarrow$ (ii):  This follows by taking cohomology.
\end{proof}

Theorem~\ref{thm:Gorenstein_DGAs} and
proposition~\ref{pro:AR_triangles_existence} combine to give:

\begin{Corollary}
\label{cor:AR_triangles_existence}
With $d = \sup \{\, i \,|\, \H^i\!R \not= 0 \,\}$,
the following conditions are equivalent.
\begin{enumerate}

  \item  Both $\Dc(R)$ and $\Dc(R^{\opp})$ have \AR\ triangles. 

  \item  There are isomorphisms of graded $\H\!R$-modules
         ${}_{\H\!R}(\dual\!\H\!R) \cong {}_{\H\!R}(\Sigma^d \H\!R)$
         and $(\dual\!\H\!R)_{\H\!R} \cong (\Sigma^d \H\!R)_{\H\!R}$.

\end{enumerate}
\end{Corollary}

\section{\AR\ triangles over a topological space}
\label{sec:AR_triangles_top}

Sections~\ref{sec:AR_triangles_top}, \ref{sec:quivers_top},
and~\ref{sec:spheres} form the topological part of this paper.  They
develop the theory of \AR\ triangles and quivers over topological
spaces, and apply the theory to spheres.  

This section proves that existence of \AR\ triangles characterizes
Poincar\'e duality spaces (theorem~\ref{thm:Poincare_duality_1}), and
gives a formula for \AR\ triangles when they exist
(proposition~\ref{pro:concrete_AR_triangle_top}).
Theorem~\ref{thm:Poincare_duality_1} is the first main result of this
paper.

\begin{Setup}
\label{set:topology}
In sections~\ref{sec:AR_triangles_top}, \ref{sec:quivers_top},
and~\ref{sec:spheres}, singular cohomology and singular cochain DGAs
are only considered with coefficients in the field $k$.  So when
$\Space$ is a topological space, $\H^{\ast}(\Space;k)$ and
$\C^{\ast}(\Space;k)$ are abbreviated to $\HXk$ and $\CXk$. Moreover,
$\D(\C^{\ast}(\Space;k))$ is abbreviated to $\D(\Space)$, and this is
combined freely with other adornments. So for instance,
$\Dc(\Space^{\opp})$ stands for $\Dc(\C^{\ast}(\Space;k)^{\opp})$.
\end{Setup}

\begin{Remark}
\label{rmk:invariance}
Recall $R$, the DGA from setup~\ref{set:blanket_R}.  If $S$ is a DGA
which is equivalent by a series of quasi-isomorphisms to $R$, then by
\cite[thm.\ III.4.2]{KrizMayAst} 
the derived categories $\D(S)$ and $\D(R)$ are equivalent.  Hence the
results of sections~\ref{sec:derived_categories},
\ref{sec:AR_triangles_DGA}, and~\ref{sec:Gorenstein_DGAs} on derived
categories apply to $S$. 

In particular, if $\Space$ is a simply connected topological space
with $\dim_k \HXk < \infty$, then $\CXk$ is equivalent by a series of
quasi-isomorphisms to a DGA
satisfying the conditions of setup~\ref{set:blanket_R}, by the methods
of \cite[proof of thm.\ 3.6]{FHT} and \cite[exam.\ 6, p.\
146]{FHTbook}. This DGA can be used as $R$ in
setup~\ref{set:blanket_R}, so the results of  
sections~\ref{sec:derived_categories}, \ref{sec:AR_triangles_DGA},
and~\ref{sec:Gorenstein_DGAs} on derived categories apply to $\CXk$. 
\end{Remark}

By this remark and lemmas~\ref{lem:Df_Krull_Schmidt}
and~\ref{lem:Dc_in_Df}, if $\Space$ is a simply connected topological
space with $\dim_k \HXk < \infty$, then $\Dc(\Space)$ and
$\Dc(\Space^{\opp})$ are categories of the sort for which
\AR\ triangles were defined in~\ref{def:AR_triangles}, so the concept
makes sense for them.

\begin{Theorem}
\label{thm:Poincare_duality_1}
Let $\Space$ be a simply connected topological space with $\dim_k
\HXk < \infty$.  Then the following conditions are
equivalent.
\begin{enumerate}

  \item  $\Space$ has Poincar\'e duality over $k$.

  \item  Both $\Dc(\Space)$ and $\Dc(\Space^{\opp})$ have \AR\ triangles. 

\end{enumerate}
\end{Theorem}

\begin{proof}
Remark~\ref{rmk:invariance} gives that
corollary~\ref{cor:AR_triangles_existence} applies to $\CXk$, the
singular cochain DGA of $\Space$ with coefficients in $k$.  For this
particular DGA, condition (ii) of
corollary~\ref{cor:AR_triangles_existence} simply says that $\Space$
has Poincar\'e duality over $k$.  So the present theorem follows.
\end{proof}

\begin{Proposition}
\label{pro:concrete_AR_triangle_top}
Let $\Space$ be a simply connected topological space with $\dim_k
\HXk < \infty$ which has Poincar\'e duality over $k$, and
write $d = \sup \{\, i \,|\, \HiXk \not= 0
\,\}$.  
\begin{enumerate}

  \item  Let $P$ be an indecomposable object of $\Dc(\Space)$.  Then
there is an \AR\ triangle in $\Dc(\Space)$,
\[
  \Sigma^{d-1}P
  \longrightarrow N
  \longrightarrow P
  \longrightarrow.
\]

  \item  The \AR\ translation of $\Dc(\Space)$ is given by 
\[
  \tau(-) = \Sigma^{d-1}(-).
\]
\end{enumerate}
\end{Proposition}

\begin{proof}
(i): Theorem~\ref{thm:Poincare_duality_1} gives that $\Dc(\Space)$ has
\AR\ triangles.  Remark~\ref{rmk:invariance} gives that
proposition~\ref{pro:AR_triangles_computation}(i) applies to
$\CXk$.  Hence there is an \AR\ triangle in $\Dc(\Space)$,
\[
  \Sigma^{-1}(\dual\!\CXk \LTensor_{\CXk} P)
  \longrightarrow N
  \longrightarrow P
  \longrightarrow.
\]
But it is easy to see from Poincar\'e duality for $X$ over $k$ that
$\dual\!\CXk$ is isomorphic to $\Sigma^d\CXk$ in the derived
category of \DGCXkbm s.  So in fact, the \AR\ triangle is the one
given in the proposition. 

\smallskip
\noindent
(ii): This is immediate from part (i); cf.\
definition~\ref{def:AR_translate}. 
\end{proof}

\section{The quiver over a topological space}
\label{sec:quivers_top}

Recall the conventions from setup~\ref{set:topology}.  When $\Space$
is a topological space, I can consider the quiver of $\Dc(\Space)$.
Moreover, when $\Space$ is simply connected with $\dim_k
\H^{\ast}(\Space) < \infty$ and with Poincar\'e duality over $k$, then
$\Dc(\Space)$ has \AR\ triangles by
theorem~\ref{thm:Poincare_duality_1} so the quiver of $\Dc(\Space)$ is
a stable translation quiver by
corollary~\ref{cor:AR_triangles_give_translation_quivers}.

\begin{Proposition}
\label{pro:quiver_weak_homotopy_invariant}
The quiver of $\Dc(\Space)$ is a weak homotopy invariant of $\Space$.

Moreover, if $\Space$ is restricted to simply connected topological
spaces with $\dim_k \HXk < \infty$ which have Poincar\'e duality over
$k$, then the quiver of $\Dc(\Space)$, viewed as a stable translation
quiver, is a weak homotopy invariant of $\Space$.
\end{Proposition}

\begin{proof}
If $\Space$ and $\Spaceprime$ have the same weak homotopy type, then
$\CXk$ and $\CXprimek$ are equivalent by a series of
quasi-isomorphisms as follows from 
\cite[thm.\ 4.15 and its proof]{FHTbook}.  Hence $\D(\Space)$ and
$\D(\Spaceprime)$ are equivalent categories by \cite[thm.\
III.4.2]{KrizMayAst}, and so the same holds for $\Dc(\Space)$ and
$\Dc(\Spaceprime)$.  This implies both parts of the proposition.
\end{proof}

\section{Spheres}
\label{sec:spheres}

Recall the conventions from setup~\ref{set:topology}.  The
$d$-dimensional sphere $S^d$ has Poincar\'e duality over any field, so
for $d \geq 2$ the category $\Dc(S^d)$ has \AR\ triangles by
theorem~\ref{thm:Poincare_duality_1}.

This section determines the \AR\ triangles in $\Dc(S^d)$ for $d \geq
2$ when $k$ has characteristic zero
(theorem~\ref{thm:AR_triangles_over_spheres}).  As a consequence
follows the determination of the quiver of $\Dc(S^d)$
(theorem~\ref{thm:quivers_over_spheres}), and it is observed that the
quiver is a sufficiently sensitive invariant to tell spheres of
different dimension apart
(corollary~\ref{cor:quivers_tell_spheres_apart}).  These are the
paper's second main results.

To determine the \AR\ triangles, I must first determine the possible
end terms, that is, the indecomposable objects of $\Dc(S^d)$.  This
requires some preparations which take up most of this section.  The
method is to set up in lemma~\ref{lem:equivalence} an equivalence of
categories between $\Dc(S^d)$ and another category whose
indecomposable objects turn out to be tractable by
lemma~\ref{lem:indecomposables_over_OmegaSd}.  Transporting these
objects through the equivalence then gives the indecomposable objects
in $\Dc(S^d)$ in proposition~\ref{pro:indecomposables_over_spheres}.

\begin{Setup}
\label{set:OmegaSd_and_E}
In this section, $d \geq 2$ is always assumed.

Let $\OmegaSd$ be the graded algebra $k[T]$ with $\deg T = -d+1$, and
view $\OmegaSd$ as a DGA over $k$ with vanishing differential.

Now $\OmegaSd/\OmegaSd^{\leq -1} \cong k$ can be viewed as a
\DGOmegaSdrm, $k_{\OmegaSd}$.  Let $\OmegaSdkres \longrightarrow
k_{\OmegaSd}$ be a $K$-projective resolution.

Let $\OmegaSdkresendo =
\Hom_{\OmegaSd^{\opp}}(\OmegaSdkres,\OmegaSdkres)$ be the endomorphism
DGA of $\OmegaSdkres$.
\end{Setup}

The point of this setup is that there is a nice connection between the
derived categories of $\OmegaSdkresendo$ and $\OmegaSd$.  Here
$\OmegaSdkresendo$ is interesting because it turns out to be
equivalent by a series of quasi-isomorphisms to $\CSdk$ when the
characteristic of $k$ is zero.  The 
algebra $\OmegaSd$ is not so interesting in itself, but is needed
because it is more computationally tractable than $\OmegaSdkresendo$
and $\CSdk$.  

The connection between the derived categories of $\OmegaSdkresendo$ and
$\OmegaSd$ can be obtained with the methods of \cite{DwyGreen}
which work because $k_{\OmegaSd}$ is a \compact\ object of
$\D(\OmegaSd^{\opp})$, as one easily checks (see also
setup~\ref{set:E}).  It takes the following form: $\OmegaSdkres$
acquires the structure $\OmegaSdkres_{\OmegaSd,\OmegaSdkresendo}$ in a
canonical way, and there are quasi-inverse equivalences of
categories,
\[
  \begin{diagram}[labelstyle=\scriptstyle]
    \tors && \pile{ \lTo^{\OmegaSdkres \LTensor_{\OmegaSdkresendo} -} \\
                    \rTo_{\RHom_{\OmegaSd^{\opp}}(\OmegaSdkres,-)} } 
    && \D(\OmegaSdkresendo), \\
  \end{diagram}
\]
where $\tors$ is a certain full triangulated subcategory of
$\D(\OmegaSd^{\opp})$ which contains $k_{\OmegaSd}$.

Since $k_{\OmegaSd}$ is in $\tors$, so is every object finitely built
from $k_{\OmegaSd}$.  It is easy to check that such objects are
exactly the ones in $\Df(\OmegaSd^{\opp})$.  Moreover, under the above
equivalences, the object $k_{\OmegaSd}$ in $\tors$ corresponds to the
object 
\[
  \RHom_{\OmegaSd^{\opp}}(\OmegaSdkres,k_{\OmegaSd}) 
  \cong \RHom_{\OmegaSd^{\opp}}(\OmegaSdkres,\OmegaSdkres) 
  \cong {}_{\OmegaSdkresendo}\OmegaSdkresendo
\]
in $\D(\OmegaSdkresendo)$, so objects finitely built from
$k_{\OmegaSd}$ correspond to objects finitely built from
${}_{\OmegaSdkresendo}\OmegaSdkresendo$.  By
lemma~\ref{lem:comp_built_from_R} these are exactly the objects of
$\Dc(\OmegaSdkresendo)$.

So the above equivalences restrict to quasi-inverse equivalences
\begin{equation}
\label{equ:equivalence}
  \begin{diagram}[labelstyle=\scriptstyle]
    \Df(\OmegaSd^{\opp}) 
    && \pile{ \lTo^{\OmegaSdkres \LTensor_{\OmegaSdkresendo} -} \\
              \rTo_{\RHom_{\OmegaSd^{\opp}}(\OmegaSdkres,-)} } 
    && \Dc(\OmegaSdkresendo). \\
  \end{diagram}
\end{equation}

To go on, it is convenient to make a specific choise of
$\OmegaSdkres$.
\begin{Setup}
\label{set:E}
Consider the morphism
\[
  \Sigma^{d-1}k[T] \longrightarrow k[T], \;\;\;
  \Sigma^{d-1}1_{k[T]} \longmapsto T
\]
of \DGrm s over $\OmegaSd = k[T]$.  Its mapping cone is easily seen to
be a minimal $K$-projective resolution of $k_{\OmegaSd}$, and from now
on I will use this mapping cone as $\OmegaSdkres$.
\end{Setup}

Observe
\begin{equation}
\label{equ:F_natural}
  \OmegaSdkres^{\natural} 
  \cong \Sigma\Sigma^{d-1}k[T]^{\natural} \amalg k[T]^{\natural}
  \cong \Sigma^d \OmegaSd^{\natural} \amalg \OmegaSd^{\natural}.
\end{equation}
Now I can prove:

\begin{Lemma}
\label{lem:E_equivalent_to_CSdk}
Suppose that $k$ has characteristic zero.  Then $\OmegaSdkresendo$ is
equivalent by a series of quasi-isomorphisms to $\CSdk$.
\end{Lemma}

\begin{proof}
The sphere $S^d$ is a so-called formal space, so since $k$ has
characteristic zero, $\CSdk$ is equivalent by a series of
quasi-isomorphisms to $\HSdk$ viewed as a DGA
with vanishing differential (see \cite[exam.\ 1, p.\ 142]{FHTbook}).
Hence it is enough to see that $\OmegaSdkresendo$ is equivalent by a
series of quasi-isomorphisms to $\HSdk$ viewed as a DGA with vanishing
differential. 

$\HSdk$ is a very simple DGA: It has a copy of $k$ in degree
zero, spanned by $1_{\HSdk}$, and another copy of $k$ in degree $d$,
spanned by some element, say $S$.

The cohomology of $\OmegaSdkresendo$ is
\[
  \H\!\OmegaSdkresendo 
  = \H(\Hom_{\OmegaSd^{\opp}}(\OmegaSdkres,\OmegaSdkres)) 
  \cong \H(\Hom_{\OmegaSd^{\opp}}(\OmegaSdkres,k_{\OmegaSd})) = (*),
\]
and as $\OmegaSdkres$ is minimal, this is
\begin{align*}
  (*) 
  & \cong \Hom_{\OmegaSd^{\opp}}(\OmegaSdkres,k_{\OmegaSd})^{\natural}
    \cong \Hom_{(\OmegaSd^{\opp})^{\natural}}
            (\OmegaSdkres^{\natural},k^{\natural}) \\
  & \stackrel{\rm (a)}{\cong}
          \Hom_{(\OmegaSd^{\opp})^{\natural}}
            (\Sigma^d \OmegaSd^{\natural} \amalg \OmegaSd^{\natural}
         ,k^{\natural})
  \cong \Sigma^{-d}k^{\natural} \oplus k^{\natural},
\end{align*}
where (a) is by equation \eqref{equ:F_natural}.  So
$\H\!\OmegaSdkresendo$ also has copies of $k$ in degrees $0$ and $d$.

Let $e$ be a cycle in $\OmegaSdkresendo^d$ whose cohomology class
spans the copy of $k$ in degree $d$ of $\H\!\OmegaSdkresendo$.  It is
now easy to check that
\[
  \HSdk \longrightarrow \OmegaSdkresendo; \; \; \; \;
  1_{\HSdk} \mapsto 1_{\OmegaSdkresendo}, \; \; S \mapsto e
\]
is a quasi-isomorphism of DGAs, proving the lemma.
\end{proof}

From lemma~\ref{lem:E_equivalent_to_CSdk} and \cite[thm.\
III.4.2]{KrizMayAst} follows that $\Dc(\OmegaSdkresendo)$ and
$\Dc(S^d)$ are equivalent.  Combining this with equation
\eqref{equ:equivalence} gives the next result.

\begin{Lemma}
\label{lem:equivalence}
Suppose that $k$ has characteristic zero.  Then there are
quasi-inverse equivalences of categories,
\[
  \begin{diagram}[labelstyle=\scriptstyle]
    \Df(\OmegaSd^{\opp}) && \pile{ \lTo \\ \rTo } 
    && \Dc(S^d).
  \end{diagram}
\]
\end{Lemma}

Let me now determine the indecomposable objects of
$\Df(\OmegaSd^{\opp})$. 

\begin{Definition}
\label{def:indecomposables_over_OmegaSd}
For each $m \geq 0$ the element $T^{m+1}$ generates a DG ideal
$(T^{m+1})$ in $k[T]$, so I can define a \DGrm\ over $\OmegaSd = k[T]$
by
\[
  \OmegaSdind_m = k[T]/(T^{m+1}).
\]
\end{Definition}

\begin{Lemma}
\label{lem:indecomposables_over_OmegaSd}
Up to isomorphism, the indecomposable objects of the category
$\Df(\OmegaSd^{\opp})$ are exactly the (positive and negative)
suspensions
\[
  \Sigma^j \OmegaSdind_m
\]
with $j$ in $\BZ$ and $m \geq 0$.
\end{Lemma}

\begin{proof}
When $K$ is a graded \OmegaSdnaturalrm, let $\delta K$
denote $K$ viewed as a \DGOmegaSdrm\ with vanishing differential.  I
claim that
\[
  K \longmapsto \delta K
\]
induces a bijective correspondence between the isomorphism clas\-ses
of $k$-finite dimensional graded indecomposable \OmegaSdnaturalrm s
and the isomorphism clas\-ses of indecomposable objects of
$\Df(\OmegaSd^{\opp})$.

For this, first note that if $M$ is a \DGOmegaSdrm, then the
cohomology $\H\!M$ is a graded right-$\H\!\OmegaSd$-module.  But
$\OmegaSd$ has vanishing differential, so $\H\!\OmegaSd$ is just
$\OmegaSd^{\natural}$, so $\H\!M$ is a graded
\OmegaSdnaturalrm.  Now in fact, I have that $M$ and
$\delta\H\!M$ are quasi-isomorphic.  This is easy to prove directly;
it is also a well known manifestation of $\OmegaSd^{\natural}$ being
graded hereditary.  (This means that any graded submodule of a graded
projective module is again graded projective.  The algebra
$\OmegaSd^{\natural}$ is graded hereditary because it is a polynomial
algebra on one generator.)  So I have $M \cong \delta\H\!M$ in
$\D(\OmegaSd)$. 

Also, if $K$ is a graded \OmegaSdnaturalrm, then I have $K \cong
\H\!\delta K$. 

Observe that this does not set up an equivalence of categories,
as the isomorphism $M \cong \delta\H\!M$ is not natural.  However,
it does show that $K \longmapsto \delta K$ induces a bijective
correspondence between the isomorphism classes of graded
\OmegaSdnaturalrm s and the isomorphism classes of
$\D(\OmegaSd^{\opp})$. 

Now, if $M$ is an indecomposable object of $\Df(\OmegaSd^{\opp})$,
then by the above I have $M \cong \delta\H\!M$ in
$\Df(\OmegaSd^{\opp})$.  If $\H\!M \cong K_1 \amalg K_2$ were 
a non-trivial decomposition, then
\[
  M \cong \delta\H\!M \cong \delta (K_1 \amalg K_2) 
  \cong \delta K_1 \amalg \delta K_2
\]
would clearly be a non-trivial decomposition in
$\Df(\OmegaSd^{\opp})$, a contradiction.  So $\H\!M$ is a $k$-finite
dimensional graded indecomposable \OmegaSdnaturalrm.

On the other hand, if $K$ is a $k$-finite dimensional graded
indecomposable \OmegaSdnaturalrm, then a similar argument shows that
$\delta K$ is an indecomposable object of $\Df(\OmegaSd^{\opp})$.

So $K \longmapsto \delta K$ also induces a bijective correspondence
between isomorphism classes of indecomposables, as claimed.

However, the finitely generated graded indecomposable
\OmegaSdnaturalrm s are exactly the (positive and
negative) suspensions of graded cyclic
\OmegaSdnaturalrm s.  This is a manifestation of
$\OmegaSd^{\natural}$ being a principal ideal domain, see \cite[p.\
9]{Spanierbook} for the ungraded case.  The $k$-finite dimensional
among these modules are
\[
  \Sigma^j(k[T]/(T^{m+1}))
\]
with $j$ in $\BZ$ and $m \geq 0$.  

By the above correspondence, up to isomorphism, the indecomposable
objects of $\Df(\OmegaSd^{\opp})$ are then
\[
  \delta \Sigma^j(k[T]/(T^{m+1}))
\]
with $j$ in $\BZ$ and $m \geq 0$.  And these are exactly the objects
$\Sigma^j \OmegaSdind_m$.
\end{proof}

Transporting the $\Sigma^j \OmegaSdind_m$'s through the equivalence of
lemma~\ref{lem:equivalence} at last gives the indecomposable objects
of $\Dc(S^d)$.

\begin{Definition}
\label{def:indecomposables_over_spheres}
Suppose that $k$ has characteristic zero.
For each $m \geq 0$ I let $\CSdkind_m$ be the object of $\Dc(S^d)$
obtained by transporting $\OmegaSdind_m$ through the equivalence of
lemma~\ref{lem:equivalence}.
\end{Definition}

\begin{Proposition}
\label{pro:indecomposables_over_spheres}
Suppose that $k$ has characteristic zero.
\begin{enumerate}

  \item  Up to isomorphism, the indecomposable objects of $\Dc(S^d)$
         are exactly the (positive and negative) suspensions 
         \[
           \Sigma^j \CSdkind_m 
         \]
         with $j$ in $\BZ$ and $m \geq 0$.

  \item  Each object of $\Dc(S^d)$ is the coproduct of uniquely
         determined indecomposable objects.

  \item  For each $m \geq 0$ the object $\CSdkind_m$ in $\Dc(S^d)$ has
\[
  \H^i\!\CSdkind_m =
  \left\{
    \begin{array}{cl}
      k & \mbox{ for $i = -m(d-1)$ and $i = d$, } \\
      0 & \mbox{ otherwise. }
    \end{array}
  \right.
\]

\end{enumerate}
\end{Proposition}

\begin{proof}
(i):  This is clear from
lemma~\ref{lem:indecomposables_over_OmegaSd} and
definition~\ref{def:indecomposables_over_spheres}. 

\smallskip
\noindent
(ii):  Remark~\ref{rmk:invariance} gives that
lemmas~\ref{lem:Df_Krull_Schmidt} and~\ref{lem:Dc_in_Df} apply to
$\CSdk$.  Hence $\Dc(S^d)$ is a Krull-Schmidt category by
\cite[3.1]{HapDerCat}, so (ii) holds.

\smallskip
\noindent
(iii):  It is easy to see that there is a distinguished triangle in
$\D(\OmegaSd^{\opp})$, 
\begin{equation}
\label{equ:OmegaSdind_triangle}
  \Sigma^{(m+1)(d-1)}\OmegaSd
  \longrightarrow \OmegaSd
  \longrightarrow \OmegaSdind_m
  \longrightarrow.
\end{equation}
It is also easy to prove
\begin{equation}
\label{equ:OmegaSd_Gorenstein}
  \H^i( \RHom_{\OmegaSd^{\opp}}(\OmegaSdkres,\OmegaSd) ) 
  \cong
  \left\{
    \begin{array}{cl}
      k & \mbox{ for $i = d$, } \\
      0 & \mbox{ otherwise. }
    \end{array}
  \right.
\end{equation}
Applying $\RHom_{\OmegaSd^{\opp}}(\OmegaSdkres,-)$ to the distinguished
triangle \eqref{equ:OmegaSdind_triangle} gives a distinguished triangle in
$\D(\OmegaSdkresendo)$, 
\[
  {\scriptstyle \Sigma^{(m+1)(d-1)}\RHom_{\OmegaSd^{\opp}}(\OmegaSdkres,\OmegaSd) }
  \longrightarrow {\scriptstyle \RHom_{\OmegaSd^{\opp}}(\OmegaSdkres,\OmegaSd) }
  \longrightarrow {\scriptstyle \RHom_{\OmegaSd^{\opp}}(\OmegaSdkres,\OmegaSdind_m) }
    \longrightarrow,
\]
and the long exact cohomology sequence and equation
\eqref{equ:OmegaSd_Gorenstein} then prove
\begin{equation}
\label{equ:intermediate_cohomology}
  \H^i( \RHom_{\OmegaSd^{\opp}}(\OmegaSdkres,\OmegaSdind_m) ) =
  \left\{
    \begin{array}{cl}
      k & \mbox{ for $i = -m(d-1)$ and $i = d$, } \\
      0 & \mbox{ otherwise. }
    \end{array}
  \right.
\end{equation}

Now, to transport $\OmegaSdind_m$ through the equivalence of
lemma~\ref{lem:equivalence} means first to transport it through the
equivalence \eqref{equ:equivalence}, secondly to transport the
resulting object through the equivalence induced by
lemma~\ref{lem:E_equivalent_to_CSdk}.  The first of these steps gives
$\RHom_{\OmegaSd^{\opp}}(\OmegaSdkres,\OmegaSdind_m)$ whose cohomology
is in equation \eqref{equ:intermediate_cohomology}.  And the second
step leaves the cohomology unchanged, viewed as a graded $k$-vector
space.  This proves the proposition's formula for $\H^i\!\CSdkind_m$.
\end{proof}

\begin{Remark}
\label{rmk:N0_equals_CSdk}
It is easy to see that $\CSdk$ itself is an indecomposable
object of $\Dc(S^d)$.  By
proposition~\ref{pro:indecomposables_over_spheres}, parts (i) and
(iii), the only possibility is
\[
  \CSdkind_0 \cong \CSdk
\]
in $\Dc(S^d)$.
\end{Remark}

Now to the second main results of this paper, which sum up the theory
in the case of spheres.  Recall from setup~\ref{set:OmegaSd_and_E} the
condition $d \geq 2$.

\begin{Theorem}
\label{thm:AR_triangles_over_spheres}
Suppose that $k$ has characteristic zero.  
\begin{enumerate}

  \item  In the category $\Dc(S^d)$, there is an \AR\ triangle
\[
  \Sigma^{d-1}\CSdkind_0 
  \longrightarrow \CSdkind_1 
  \longrightarrow \CSdkind_0 
  \longrightarrow
\]
and an \AR\ triangle
\[
  \Sigma^{d-1}\CSdkind_{n} 
  \longrightarrow \Sigma^{d-1}\CSdkind_{n-1} \oplus \CSdkind_{n+1}
  \longrightarrow \CSdkind_{n}
  \longrightarrow
\]
for each $n$ with $n \geq 1$, where the $\CSdkind$'s are the
indecomposable objects from
definition~\ref{def:indecomposables_over_spheres}.  Each
\AR\ triangle is a (positive or negative) suspension of one of these.

  \item  The \AR\ translation of $\Dc(S^d)$ is given by
\[
  \tau(-) = \Sigma^{d-1}(-).
\]
\end{enumerate}
\end{Theorem}

\begin{proof}
(i): By \cite[prop.\ 3.5(i)]{HapDerCat}, \AR\ triangles are determined
up to isomorphism by their right hand end terms.  The right hand end
terms are indecomposable objects by definition, so in the present case
have the form $\Sigma^j \CSdkind_m$ with $j$ in $\BZ$ and $m \geq 0$
by proposition~\ref{pro:indecomposables_over_spheres}(i).  So to prove
part (i) of the theorem, it is clearly enough to see that the \AR\
triangles with right hand end terms $\CSdkind_{m}$ for $m
\geq 0$ are as claimed.

By proposition~\ref{pro:concrete_AR_triangle_top}(i), the left hand end
terms of the \AR\ triangles are as claimed in the theorem, so let me
consider the middle terms.  First the \AR\ triangle ending in
$\CSdkind_0$,
\begin{equation}
\label{equ:N0_triangle}
  \Sigma^{d-1}\CSdkind_0 
  \longrightarrow N
  \longrightarrow \CSdkind_0 
  \stackrel{\pi}{\longrightarrow}.
\end{equation}
By definition~\ref{def:AR_triangles} the morphism $\pi$ is non-zero.
But by remark~\ref{rmk:N0_equals_CSdk} this morphism is $\CSdk
\stackrel{\pi}{\longrightarrow} \Sigma^d \CSdk$.
This makes it easy to compute the long exact cohomology sequence of
\eqref{equ:N0_triangle} and get
\[
  \H^i\!N
  \cong
  \left\{
    \begin{array}{cl}
      k & \mbox{ for $i=-(d-1)$ and $i = d$, } \\
      0 & \mbox{ otherwise. }
    \end{array}
  \right.
\]
But $N$ is the coproduct of uniquely determined indecomposable objects
of $\Dc(S^d)$ by
proposition~\ref{pro:indecomposables_over_spheres}(ii), and
by~\ref{pro:indecomposables_over_spheres}, parts (i) and (iii), the
only possibility is $N \cong \CSdkind_1$.

Next the \AR\ triangle ending in $\CSdkind_{n}$,
\begin{equation}
\label{equ:Nell_triangle}
  \Sigma^{d-1}\CSdkind_{n}
  \longrightarrow N
  \longrightarrow \CSdkind_{n}
  \stackrel{\pi}{\longrightarrow},
\end{equation}
with $n \geq 1$.  There can be no retract $\Sigma^j \CSdk
\longrightarrow \CSdkind_{n}$, for else
$\CSdkind_{n}$ would be a direct summand in the indecomposable object
$\Sigma^j \CSdk \cong \Sigma^j \CSdkind_0$.  Hence each morphism
$\Sigma^j \CSdk \stackrel{\gamma}{\longrightarrow} \CSdkind_{n}$ has
$\pi\gamma = 0$.  But this shows $\H\!\pi = 0$, so the long exact
cohomology sequence of \eqref{equ:Nell_triangle} splits into short
exact sequences.  So using
proposition~\ref{pro:indecomposables_over_spheres}(iii), the cohomology
of $N$ can be read off as
\[
  \H^i\!N
  \cong
  \left\{
    \begin{array}{cl}
      k & \mbox{ for $i$ in $\{ -(n+1)(d-1), -n(d-1), 1, d \}$, } \\
      0 & \mbox{ otherwise. }
    \end{array}
  \right.
\]
Proposition~\ref{pro:indecomposables_over_spheres}(ii) says that $N$
is the coproduct of uniquely determined indecomposable objects of
$\Dc(S^d)$.  Comparing the cohomology of $N$ with the cohomology of
the indecomposable objects, obtained
from~\ref{pro:indecomposables_over_spheres}, parts (i) and (iii),
leaves only two possibilities: $N$ is either
$\Sigma^{d-1}\CSdkind_{n-1} \oplus
\CSdkind_{n+1}$ or $\Sigma^{d-1}\CSdkind_{n} \oplus
\CSdkind_{n}$. 

However, let me suppose by induction that the \AR\ triangle ending in
$\CSdkind_{n-1}$ is as claimed in the theorem, hence has a
summand $\CSdkind_{n}$ in its middle term.  By 
lemma~\ref{lem:AR_triangles_determine_quivers}, (iii) $\Rightarrow$ (ii),
this implies that there is an irreducible morphism
$\CSdkind_{n} \longrightarrow \CSdkind_{n-1}$.  Hence there is
an irreducible morphism $\Sigma^{d-1}\CSdkind_{n} \longrightarrow
\Sigma^{d-1}\CSdkind_{n-1}$, and by
lemma~\ref{lem:AR_triangles_determine_quivers}, (i) $\Rightarrow$ (iii),
this implies that $\Sigma^{d-1}\CSdkind_{n-1}$
is a direct summand of $N$.  So $N$ must be
$\Sigma^{d-1}\CSdkind_{n-1} \oplus \CSdkind_{n+1}$, proving the
theorem.

(ii):  This is immediate from part (i), or from
proposition~\ref{pro:concrete_AR_triangle_top}(ii).
\end{proof}

If a category has \AR\ triangles, then knowledge of the \AR\ triangles
gives full knowledge of the quiver of the category by
lemma~\ref{lem:AR_triangles_determine_quivers}.  Also in this case,
the quiver is a stable translation quiver with translation induced by
the \AR\ translation of the category, by
corollary~\ref{cor:AR_triangles_give_translation_quivers}.  Applying
this to the data from theorem~\ref{thm:AR_triangles_over_spheres}
gives the following.

\begin{Theorem}
\label{thm:quivers_over_spheres}
Suppose that $k$ has characteristic zero.  Then the quiver of the
category $\Dc(S^d)$ consists of $d-1$ components, each isomorphic to
$\BZ \BA_{\hspace{0.2mm}\infty}$.  The component containing
$\CSdkind_0 \cong \CSdk$ is
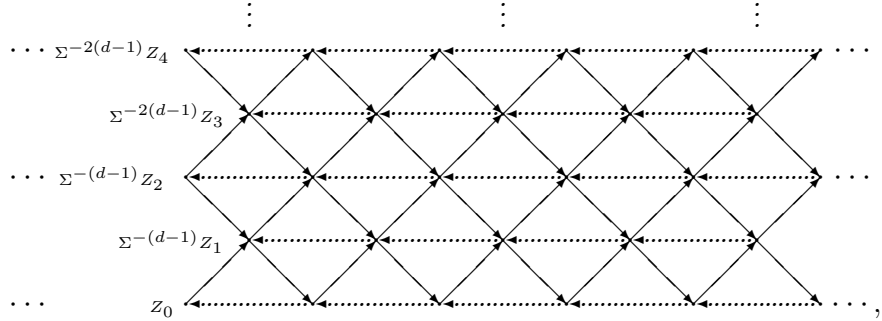
\[
  \begin{diagram}[height=1em,width=1em,abut]
    &&&&&&& \vdots &&&&&&&& \vdots &&&&&&&& \vdots &&& \\
    \cdots & \lefteqn{\scriptscriptstyle \Sigma^{-2(d-1)}\CSdkind_4} \; &&&& \cdot && \lDotsto && \cdot && \lDotsto && \cdot && \lDotsto && \cdot && \lDotsto && \cdot && \lDotsto && \cdot & \cdots \\
    &&&&&& \rdTo && \ruTo && \rdTo && \ruTo && \rdTo && \ruTo && \rdTo && \ruTo && \rdTo && \ruTo && \\
    &&& \lefteqn{\scriptscriptstyle \Sigma^{-2(d-1)}\CSdkind_3} \;\;\; &&&& \cdot && \lDotsto && \cdot && \lDotsto && \cdot && \lDotsto && \cdot && \lDotsto && \cdot &&& \\
    &&&&&& \ruTo && \rdTo && \ruTo && \rdTo && \ruTo && \rdTo && \ruTo && \rdTo && \ruTo && \rdTo && \\
    \cdots & \lefteqn{\scriptscriptstyle \Sigma^{-(d-1)}\CSdkind_2} &&&& \cdot && \lDotsto && \cdot && \lDotsto && \cdot && \lDotsto && \cdot && \lDotsto && \cdot && \lDotsto && \cdot & \cdots \\
    &&&&&& \rdTo && \ruTo && \rdTo && \ruTo && \rdTo && \ruTo && \rdTo && \ruTo && \rdTo && \ruTo && \\
    &&& \lefteqn{\scriptscriptstyle \Sigma^{-(d-1)}\CSdkind_1} \; &&&& \cdot && \lDotsto && \cdot && \lDotsto && \cdot && \lDotsto && \cdot && \lDotsto && \cdot &&& \\
    &&&&&& \ruTo && \rdTo && \ruTo && \rdTo && \ruTo && \rdTo && \ruTo && \rdTo && \ruTo && \rdTo && \\
    \cdots &&&& \lefteqn{\scriptscriptstyle \CSdkind_0} \; & \cdot && \lDotsto && \cdot && \lDotsto && \cdot && \lDotsto && \cdot && \lDotsto && \cdot && \lDotsto && \cdot & \cdots \lefteqn{,} \\
  \end{diagram}
\]
where the unbroken arrows are the arrows of the quiver and the
dotted arrows indicate the action of the translation induced by the
\AR\ translation of $\Dc(S^d)$.
\end{Theorem}

Finally, the following corollary is clear from
theorem~\ref{thm:quivers_over_spheres}. 

\begin{Corollary}
\label{cor:quivers_tell_spheres_apart}
Suppose that $k$ has characteristic zero.  Then the quiver of
$\Dc(S^d)$ is a sufficiently sensitive invariant to tell different
$S^d$'s apart.
\end{Corollary}

\medskip
\noindent
{\bf Acknowledgement. }
Theorem~\ref{thm:Poincare_duality_1} is inspired by Happel's
result \cite[thm. 3.4]{HapGor}, which considers a finite dimensional
algebra $\Lambda$ and says roughly that $\Dc(\Lambda)$ has
\AR\ triangles if and only if $\Lambda$ is Gorenstein. This is related
to theorem~\ref{thm:Poincare_duality_1} because the differential
graded analogue of the Gorenstein property is Poincar\'e duality (see
section~\ref{sec:Gorenstein_DGAs}).  I thank Henning Krause for
directing my attention to \cite{HapGor}.

The diagrams were typeset with Paul Taylor's {\tt diagrams.tex}.

\end{document}